\documentclass[12pt]{amsart}
\usepackage{a4wide}
\usepackage{amssymb}

\usepackage{color}

\newtheorem{teo}{Theorem}
\newtheorem{defi}[teo]{Definition}
\newtheorem{eje}[teo]{Example}
\newtheorem{lem}[teo]{Lemma}
\newtheorem{pro}[teo]{Proposition}
\newtheorem{coro}[teo]{Corolary}
\newtheorem{obs}[teo]{Observation}
\newtheorem{fact}[teo]{Fact}
\newtheorem{rem}[teo]{Remark}
\newtheorem{proo}[teo]{Proof}

\def\bteo{\begin{teo}}
\def\eteo{\end{teo}} 
 \def\bdefi{\begin{defi}}
\def\edefi{\end{defi}} 
\def\beje{\begin{eje}}
\def\eeje{\end{eje}} 
\def\blem{\begin{lem}}
\def\elem{\end{lem}} 
\def\bpro{\begin{pro}}
\def\epro{\end{pro}} 
\def\bproo{\begin{proo}}
\def\eproo{\end{proo}}
\def\bcoro{\begin{coro}}
\def\ecoro{\end{coro}}     
\def\bob{\begin{obs}}
\def\eob{\end{obs}}
\def\bfact{\begin{fact}}
\def\efact{\end{fact}}
\def\brem{\begin{rem}}
\def\erem{\end{rem}}

\def\square{\ifmmode\sqr\else{$\sqr$}\fi}
\def\sqr{\vcenter{
         \hrule height.1mm
         \hbox{\vrule width.1mm height2.2mm\kern2.18mm\vrule width.1mm}
         \hrule height.1mm}} 

\def\reff#1{(\ref{#1})}

\def\N{{\mathbb N}}
\def\Z{{\mathbb Z}}
\def\R{{\mathbb R}}
\def\E{{\mathbf E}}
\def\P{{\mathbf P}}

\def\cR{{\mathcal R}}

\def\cT{{\mathcal T}}
\def\cU{{\mathcal U}}
\def\cK{{\mathcal K}}
\def\cG{{\mathcal G}}
\def\cB{{\mathcal B}}

\def\cB{{\mathcal B}}

\def\cP{{\mathcal P}}

\def\one{{\mathbf 1}}

\def\cE{{\mathcal E}}
\def\cX{{\mathcal X}}
\def\cC{{\mathcal C}}
\def\hR{{\hat R}}

\def\np{{\mathbf{p}}}
\def\nn{{\mathbf{n}}}

\title{Absence of Percolation in the Bernoulli Boolean Model}
\author{Cristian F. Coletti and Sebastian P. Grynberg}

\address{
\newline
UFABC - Centro de Matem\'atica, Computa\c{c}\~ao e Cogni\c{c}\~ao.
\newline
Avenida dos Estados, 5001, Santo Andr\'e - S\~ao Paulo, Brasil
\newline
e-mail:  \rm \texttt{cristian.coletti@ufabc.edu.br}
}

\address{
\newline
Universidad de Buenos Aires,Departamento de Matem\'aticas, Facultad de Ingenier\'ia
\newline
Av. Paseo Col\'on 850 - C1063ACV - Buenos Aires - Argentina
\newline
e-mail:  \rm \texttt{sebgryn@fi.uba.ar}
}


\begin{document}
\maketitle
\begin{abstract}
We consider the Bernoulli Boolean discrete percolation model on the d-dimensional integer lattice. We study sufficient conditions on the distribution of the radii of balls
placed at the points of a Bernoulli point process for the absence of percolation, provided that the intensity of the underlying point process is small enough. We also study a Harris
graphical procedure to construct, forward in time, particle systems with interactions of infinite range under the assumption that the corresponding generator admits a Kalikow-type
decomposition. We do so by using the subcriticality of the boolean model of discrete percolation.

\end{abstract}

\section*{Introduction.}
We study the Bernoulli Boolean discrete percolation model on the $d$-dimensional lattice $\Z^d$. This is a discrete percolation model which can be informally described
as follows. Consider a Bernoulli point process $\cX$ with retention parameters $0<p_x<1$, $x\in\Z^d$, on the $d$-dimensional lattice $\Z^d$. This means that each site
$x\in\Z^d$ is {\it{present}} or {\it{absent}} in $\cX$ with probability $p_x$ or $1-p_x$, respectively and independently of anything else. Each point of $\cX$ is the
center of a ball of random radius in the metric induced by the $L_1$ norm. The random radii $R_x$, $x\in\Z^d$, are independent and independent of $\cX$. We consider the 
\emph{occupied region} which is defined as the subset of $\Z^d$ obtained by taking the union of all random balls centred at the points of $\cX$.

This model is the discrete counterpart of the Poisson Boolean model of continuum percolation. In the  Poisson Boolean model a ball of random radius is centred at each point of a
homogeneous Poisson point process with density $\lambda$ on $\R^d$. The corresponding radii form an independent and identically distributed collection of non-negative
random variables which are also independent of the point process. Denote by $\cB$ the union of these balls and by $\cC$ the connected component of $\cB$ containing the origin. Let $R$ be one of the random radii and denote by $\bf{P}$ the law governing
the continuous boolean model. Also, denote by $\bf{E}$ the corresponding expectation operator. In \cite{Hall}, Hall proved that
for values of $\lambda$ small enough, $\cC$ is almost surely bounded provided that $\E[R^{2d-1}]$ is finite. In \cite{Meester_and_Roy}, Meester and Roy proved that if $d \geq 2$ and
$\E[R^{d}]$ is finite, then the expected number of balls in the occupied component  which contains the origin is finite whenever $\lambda$ is small enough if, and only if,
$\E[R^{2d}]$ is finite. Also, they proved that if $\E[R^{2d-1}]$ is finite then $\P(\mbox{number of balls in any occupied component is finite})=1$ provided that $\lambda$ is small
enough.  In \cite{Gouere}, Gouere showed that the set $\mathcal{C}$ is almost surely bounded for small enough $\lambda$ if and only if $\E[R^d]$ is finite. 


In this paper we prove that if $p_x=p \in (0,1)$ for all $x$ and the random radii $(R_x, x \in \mathbb{Z}^d)$ are i.i.d. random variables with finite $d$-moment, then the connected
components arising in the discrete Boolean model are almost surely finite for sufficientlu small values of $p$. We also prove that such behavior does not occur if the
random radii have infinite $d$-moment. Then, using a coupling argument, we extend the result about subcriticality to the case where the values of $p_x$ are not constant and
the random radii are independent but not necessarily identically distributed. Then we use the result above about subcriticality to provide a graphical construction method for
interacting particle systems with interactions of infinite range. In order to prove this result we assume that the generator of the particle system admits a Kalikow-type decomposition. 
Recently, this type of decomposition has been explored by Galves et al. in the
context of perfect simulation of interacting particle systems with interactions of infinite range. More precisely, in \cite{Galves_Garcia_Locherbach} the authors exhibit a
sufficient condition under which a Kalikow-type decomposition holds for the transition rates of interacting particle systems with interactions of infinite range. Namely, if the
transition rates satisfy a continuity condition then the referred decomposition holds. For further details on Kalikow-type decompositions see \cite{kalikow} and \cite{bramson-kalikow}. 

This paper is organized as follows: In section \ref{booleanpercolation} we describe the discrete boolean percolation model and state the main result of this work which
says about the absence of percolation on the model described above. This result is proved in section \ref{mainp} following ideas for the continuous boolean percolation
model studied in \cite{Gouere}. 
In section \ref{PoTP} we extend the result in \cite{harris} on the graphical
construction of interacting particle system with finite-range interaction to the case of interactions of infinite range, using the results in section \ref{booleanpercolation}
under mild assumption on the decay of the range of interaction. 

\section{Definitions, notation and main results} 
\label{booleanpercolation}
Throughout this paper $\N_0$ will denote the set of non-negative integer numbers. We write $\| \;\|$ for the $L_1$ norm on $\Z^d$ and $|A|$ for the cardinal number
of any set $A\subset\Z^d$. Also, $B(x,r)=\{y\in\Z^d:\, \|y-x\|\leq r\}$ denotes the (close) ball of radius $r$ centred at $x$ and $S_r=\{x\in\Z^d:\,\|x\|=r\}$ denotes the sphere of
radius $r$. For any set $A\subset\Z^d$, $A^c$ stands for the complement of $A$.

If $F$ denotes a cumulative distribution function, let $F^{-1}$ be the generalized inverse of $F$ defined by $F^{-1}(u)=\inf\{r\in\R: F(r)\geq u\}$ where $u\in[0,1]$. 
If $X$ and $Y$ are two stochastic elements equally distributed, we write $X \stackrel{D}{=}Y$.


A \emph{Bernoulli point process on} $\Z^d$ with retention parameters $\np=(p_x:x\in\Z^d)$, where $0<p_x<1$ for all $x\in\Z^d$, is a family of independent $\{0,1\}$-valued random variables
$\cX=(X_x:\,x\in\Z^d)$ such that $p_x$ is the probability of the event $\{X_x=1\}$. Identify the family of random variables $\cX$ with the random subset $\cP$ of $\Z^d$ defined by
$\cP=\{x\in\Z^d:\, X_x = 1\}$ whose distribution is a product measure whose marginals at each site $x$ are Bernoulli distribution of parameter $p_x$.

By a \emph{Bernoulli marked point process} on $\Z^d$ we mean a pair $(\cX,\cR)$ formed by a Bernoulli point process $\cX$ on $\Z^d$ and a family of independent $\mathbb{N}_0$-valued random variables
$\cR=(R_x: x\in\Z^d)$ called marks. We assume that these marks are independent of the point process $\cX$.    

Let $(\cX,\cR)$ be a Bernoulli marked point process on $\Z^d$. Let $p_x$ be the retention parameter of the random variable $X_x$ and let $\nu_x$ be the probability function of the random variable $R_x$. If there exists a value $p\in(0,1)$ and a probability function $\nu$ on $\N_0$ such that $p_x=p$ and $\nu_x=\nu$ for every $x\in\Z^d$ we say that the marked point process $(\cX,\cR)$ is \emph{spatially homogeneous with retention parameter $p$ and marks distributed according to $\nu$}. 

Let $(\cX,\cR)$ be a Bernoulli marked point process on $\Z^d$ with retention parameters $\np=(p_x:x\in\Z^d)$ and marks distributed according to a family of probability functions
$\mathbf{n}=(\nu_x:x\in\Z^d)$. We denote by $\P_{\np,\,\nn}$ and $\E_{\np,\,\nn}$ respectively the probability measure and the expectation operator induced by $(\cX,\cR)$. 
If $(\cX,\cR)$ is spatially homogeneous with retention parameter $p$ and marks distributed according to the probability function $\nu$, we denote by $\P_{p,\,\nu}$ and $\E_{p,\,\nu}$
respectively the probability measure and the expectation induced by $(\cX,\cR)$.

Let $(\cX,\cR)$ and $(\cX',\cR')$ be two marked point process on $\Z^d$ defined on the same probability space. If  
\[X_x\leq X'_x\qquad \mbox{ and }\qquad R_x\leq R'_x,\qquad x\in\Z^d,\]  
we say that $(\cX,\cR)$ is \emph{dominated} by $(\cX',\cR')$ and we denote this by $(\cX,\cR)\preceq (\cX',\cR')$.

\subsection*{Random Graphs and Percolation.} Let $(\cX,\cR)$ be a Bernoulli marked point process on $\Z^d$. Then we define an associated random graph $\cG(\cX,\cR)=(\Z^d,\cE)$ as the
undirected random graph with vertex set $\Z^d$ and edge set $\cE$ defined by the condition $\{x,y\}\in\cE$ if, and only if, $X_x=1$ and $y\in B(x,R_x)$ or $X_{y}=1$ and $x\in B(y,R_y)$. 

A \emph{path} on $\cG(\cX,\cR)$ is a sequence of distinct vertex $x_0, x_1,\dots, x_n$ with $x_{i-1} \neq x_i$ such that $\{x_{i-1},x_i\}\in\cE$, $i=1,\dots, n$.

A set of vertex $C\subset\Z^d$ is connected if, for all pair of distinct vertex $x$ and $y$ in $C$, there exists a path on $\cG(\cX,\cR)$ using vertices only from $C$, starting at $x$
and ending at $y$. The connected components of the graph $\cG(\cX,\cR)=(\Z^d,\cE)$ are its maximal connected subgraphs.   

The cluster $C(x)$ of vertex $x$ is the connected component of the graph $\cG(\cX,\cR)$ containing $x$. Define the Percolation as follows:
\begin{eqnarray}
[\mbox{Percolation}]:=\bigcup_{x\in\Z^d}\left\{|C(x)|=\infty\right\}.  
\end{eqnarray}
\noindent {\bf{Phase transition}}. Consider the Bernoulli Boolean discrete percolation model introduced above. Then replace the random radii in this model by the deterministic radius $0$. What
we get is the independent site percolation model. It is well known (see Grimmett \cite{Grimmett}, page 25) that the critical parameter for this last model is a positive number
$p_c^{\mbox{site}}(\Z^d)<1$. Then, a coupling argument shows that for any $p>p_c^{\mbox{site}}(\Z^d)$ there is percolation for the discrete Boolean model. Thus we focus our attention
in the subcritical regime.

\medskip

Now we state the main result of this work.

\bteo
Let $(\cX,\cR)$ be a spatially homogeneous marked point process on $\Z^d$ with retention parameter $p$ and marks distributed according to a probability function $\nu$. If $\sum_{r\geq 1}r^d\nu(r)<\infty$, then there exists $p_0>0$ such that $\P_{p,\,\nu}($Percolation$)=0$ for all  $p\leq p_0$. 
\label{coupI}
\eteo

Indeed, a similar result holds if we only assume that the values of $p_x$ are uniformly bounded and the family of random radii are independent, but not identically distributed.

\bteo
\label{Rind}
Let $(\cX,\cR)$ be a marked point process on $\Z^d$ with retention parameters $\np=(p_x:x\in\Z^d)$, and random marks $R_x$, $x\in\Z^d$ distributed according to a family of probability functions $\mathbf{n}=(\nu_x:x\in\Z^d)$. If
\begin{equation}
\label{StD}
\lim_{r\to\infty}\inf_{x\in\Z^d}\P_{\np,\,\nn}(R_x \leq r)=1
\end{equation}
and 
\begin{equation}
\label{coup} 
\sum_{r}r^d\left(\inf_{x\in\Z^d}\P_{\np,\,\nn}(R_x\leq r)-\inf_{x\in\Z^d}\P_{\np,\,\nn}(R_x\leq r-1)\right)<\infty,
\end{equation}
then there exists $p_0>0$ such that $\P_{\np,\,\nn}($Percolation$)=0$ for any family of retention parameters $\np=(p_x:x\in\Z^d)$ such that $\sup_{x\in\Z^d}p_x\leq p_0$.   
\eteo

\brem

We claim that hypothesis (\ref{StD}) in Theorem \ref{Rind} above is equivalent to assume the existence of a random variable $R$ such that each random variable in
$\mathcal{R}$ is stochastically dominated by $R$. Indeed, let $U$ be a uniform random variable on $[0,1]$. Then, define $R$ as follows:
\begin{equation}
\label{VDM}
R=\sum_{r\geq 1}r\cdot\one\left\{\inf_{x\in\Z^d}\P_{\np,\,\nn}\left(R_{x}\leq r-1\right) <U \leq \inf_{x\in\Z^d} \P_{\np,\,\nn}\left(R_{x} \leq r \right) \right\}.
\end{equation}
 We readily verify that $R$ is a random variable and that each random variable in $\mathcal{R}$ is stochastically dominated by $R$. By hypothesis (\ref{coup}) we  have $\E[R^d]<\infty$. Now, using this stochastic domination we may construct, in a common probability space, two marked point processes. For that purpose, let  $(U_x:\,x\in\Z^d)$ be a family of i.i.d. random variables, each one uniformly distributed on $[0,1]$. Then, set $\hat{R}_x=F_{R_x}^{-1}(U_x)$ and  $\hat{R}^x=F_{R}^{-1}(U_x)$, where $F_{R_x}$ and $F_{R}$ are the cumulative distribution function of $R_x$ and $R$ respectively. It follows that  $\hat{R}_x\stackrel{D}{=}R_x$, $\hat{R}^x\stackrel{D}{=}R$  and $\hat{R}_x\leq \hat{R}^x$. Since $\hat{R}_x\leq \hat{R}^x$, we get that $B(x,\hat{R}_x)\subset B(x,\hat{R}^x)$. Then,  $\left(\mathcal{X},\mathcal{R}_1\right)\preceq \left(\mathcal{X},\mathcal{R}_2\right)$, where $\cR_1=(\hat{R}_x:\, x\in\Z^d)$ and $\cR_2=(\hat{R}^x:, x\in\Z^d)$. Since $(\cX,\cR)\stackrel{D}{=}(\cX,\cR_1)$ and $\E[R^d]<\infty$, Theorem \ref{Rind} is a simple consequence of Theorem \ref{coupI}. 

\erem

\beje

Observe that condition (\ref{coup}) in Theorem \ref{Rind} turns out to be slightly stronger than requiring $\sup_{x\in\Z^d}\E_{\np,\,\nn}[R_x^d]<\infty$. Note that if the
random radii $R_x, x\in\Z^d$ are i.i.d random variables, then condition (\ref{coup}) in Theorem \ref{Rind} becomes $\E_{\np,\,\nn}[R^d]<\infty$, where $R$ is a random
variable distributed as $R_x$ for some $x\in\Z^d$. The example below shows that it is possible to construct a sequence of random variables $(R_n)_{n\in\N}$ In a common
probability space such that
\begin{eqnarray}
\sup_{n\in\N}\E[R_n]<\infty \mbox{ with } \E[R]=\infty,
\end{eqnarray}
\noindent where $R$ is a random variable such that $\P(R\leq r)=\inf_{n\in\N}\P(R_n\leq r)$ and $\E$ is the corresponding expectation operator.

Let $(R_n)_{n\geq 2}$ be a sequence of random variables with distribution function
\begin{eqnarray*}
F_{n}(x)
&=&
\left(1-\frac{3}{4n}\right)\one\{0\leq x<1\}\\
&+&
\left(\frac{1}{4n(n-1)}(x-1)+1-\frac{3}{4n}\right)\one\{1\leq x<n\}\\
&+&
\left(1-\frac{1}{2n}\right)\one\{n\leq x<n+1\}\\
&+&
\one\{n+1\leq x\}.
\end{eqnarray*} 

We readily check that

\begin{eqnarray*}
\sup_{n\geq 2}\E[R_n]<\infty.
\end{eqnarray*}

The distribution function $F(x)=\inf_{n\in\N}F_n(x)$ is given by

\begin{eqnarray*}
F(x)=\sum_{n\geq 2}\left(1-\frac{1}{2n}\right)\one\{n\leq x<n+1\}.
\end{eqnarray*}

Finally, note that if $R$ is a random variable with distribution function $F$ as above, then

\[
\E[R]\leq \sum_{n \geq 2}\frac{1}{2n}=+\infty.
\]

This example shows that, with our techniques, the hypothesis in Theorem \ref{Rind} can not be weakened.
\eeje

\subsection*{Complete Coverage}

We complement the result of Theorem \ref{coupI} by establishing a sufficient condition for complete coverage of the space $\Z^d$. For
any $A\subset\Z^d$, define $\Lambda(A)=\bigcup_{x\in A\cap\cP}B(x,R_x)$. 

\bteo
\label{Riidinf} 
Let $(\cX,\cR)$ be a spatially homogeneous marked point process on $\Z^d$ with retention parameter $p$ and marks distributed
according to the probability function $\nu$. If $\sum_{r\geq 1}r^d\nu(r)=\infty$,  then for any $p\in(0,1]$ , $\Lambda
(\mathbb{Z}^d)=\Z^d$.
\eteo

\subsection{Particle systems with interactions of infinite range}
Let $S$ be a finite (or countable) set and let $S^{\Z^d}$ be the set of mappings $\sigma:\Z^d\to S$. Give $S$ the discrete topology and $S^{\Z^d}$ the product topology. The measurable
sets of $S^{\Z^d}$ are the Borel sets. The elements of $S$ are called \emph{spins} or \emph{particles}. $S^{\Z^d}$ is called the configuration space and its elements are in general
written as $\sigma$, $\eta$, $\xi\dots$. For each $x\in\Z^d$, $\sigma(x)$ denotes the spin value of configuration $\sigma$ at site $x$. For each $A\subset\Z^d$, $\sigma(A)\in S^A$
denotes the restriction of configuration $\sigma$ to $A$.

A particle system with interactions of infinite range is a Markov process on $S^{\Z^d}$ whose generator is defined on cylinder functions by
\begin{eqnarray}
\label{IG}
Lf(\sigma)
&=&\sum_{x\in\Z^d}\sum_{s\in S}c_x(s,\sigma)\left[f(\sigma_{x,s}) - f(\sigma)\right],
\end{eqnarray}
where $\sigma_{x,s} \in S^{\mathbb{Z}^d}$ is defined by $\sigma_{x,s}(x)=s$, $\sigma_{x,s}(y)=\sigma(y)$ if $y\neq x$. Here $c_{x}(s,\sigma)>0$ is the intensity for a jump
$\sigma\to\sigma_{x,s}$ and it depends on $x$ and the whole spin configuration $\sigma$.  

\paragraph{Kalikow-type decomposition.} We assume that the following \emph{Kalikow-type decomposition} for the jump intensities holds:
\begin{eqnarray}
\label{rates}
c_x(s,\sigma)=M_xp_{x}(s|\sigma),
\end{eqnarray}
where $M_x>0$ and
\begin{eqnarray}
\label{CP}
p_{x}(s|\sigma)=\sum_{r\geq 0}\nu_x(r)p_x^{[r]}(s|\sigma). 
\end{eqnarray}   
Here, $\nu_x(\cdot)$ is a probability function  on $\N_0$ and $p_x^{[r]}(\cdot|\sigma)$ is a probability function on $S$ which depends on $\sigma$ only through $\{\sigma(y):\, y\in B(x,r)\}$. For further details on this kind of decomposition see \cite{kalikow} and \cite{bramson-kalikow}.

Here and for the rest of the paper we will assume that
\begin{equation} 
\label{supinf}
0<M_*:=\inf_{x\in\Z^d}M_x\leq M^*:=\sup_{x\in\Z^d}M_x<\infty
\end{equation}
and 
\begin{equation} 
\label{momentcondition}
\sum_{r\geq 1}r^d\sup_{x\in\Z^d}\nu_x(r)<\infty.
\end{equation}
Now we can state the result about existence of interacting particle systems with interactions of infinite range. 

\bteo
\label{Harrisconstruction}
Let $\{c_x(s,\sigma):\, x\in\Z^d,\, s\in S,\, \sigma\in S^{\Z^d}\}$ be a family of jump intensities satisfying the Kalikow-type decomposition described in \reff{rates} and \reff{CP}. Let assumptions \reff{supinf} and \reff{momentcondition} hold. Then, for each initial spin configuration $\eta$, there exists an almost surely unique interacting particle system $(\sigma_t^{\eta})_{t\geq 0}$ with generator
\begin{eqnarray}
\label{2.12}
Lf(\sigma)&=&\sum_{x\in\Z^d}\sum_{s\in S}\sum_{r \geq 0}M_x\nu_x(r)p_x^{[r]}(s|\sigma)\left[f(\sigma_{x,s}) - f(\sigma)\right].
\end{eqnarray}
\eteo

In the rest of this section we present the sketch of the proof of Theorem \ref{Harrisconstruction} based on a extension, to the case of interactions with infinite range, of ideas developed by Harris for the case of finite range interactions. See \cite{harris} for further details.    

\subsection*{Harris graphical construction.} The probability space where the Markov processes $(\sigma^{\eta}_t: \,t\geq 0)$ will be constructed is the space generated by
a family $(\cT,\cK,\cU)=\{(\cT_x, \cK_x, \cU_x):\,x\in\Z^d\}$ of mutually independent marked Poisson point processes on the time line $[0,\infty)$. For each $x\in\Z^d$,
the Poisson process $\cT_x=(T_{x,n}:\,n\in\N)$ is homogeneous with rate $M_x$, $\cK_x=(K_{x,n}:\,n\in\N)$ is a sequence of i.i.d. random variables with common law
$\nu_x$ on $\mathbb{N}_0$ and $\cU_x=(U_{x,n}:\,n\in\N)$ is a sequence of i.i.d. uniform random variables on $[0,1]$. Moreover, for each
$x \in \mathbb{Z}^d , \cT_x, \cK_x$ and $\cU_x$ are mutually independent.

For each $\eta\in S^{\Z^d}$, we construct a process $(\sigma^{\eta}_t:\,t\geq 0)$ with generator \reff{2.12} and initial condition $\eta$ at time $0$ as a function of the
family $(\cT,\cK, \cU)$. Roughly speaking, the process $(\sigma^{\eta}_t:\,t\geq 0)$ is constructed as follows. Initially, $\sigma^{\eta}_0:=\eta$. Then, at the time epoch
$t\in\cT_x$, the spin value at site $x$ is updated in the following way: if $t=T_{x,n}$, then sample the range of interaction using the random variable 
$K_{x,n}\in\cK_x$. If $K_{x,n}=r$, then the spin value at site $x$ is updated by a random variable $W_{x}(\sigma^{\eta}_{t-})$ with law
$p_{x}^{[r]}(\cdot|\sigma^{\eta}_{t-})$:
\begin{eqnarray}
\label{GSa}
\sigma^{\eta}_t=\sigma^{\eta}_{x,W_{x}(\sigma^{\eta}_{t-})}.
\end{eqnarray}
The random variable $W_{x}(\sigma^{\eta}_{t-})$ is constructed as a function of the uniform random variable $U_{x,n}\in\cU_x$.

Since there are infinitely many Poisson processes the main difficulty in the construction described above is that in general there will be infinitely many jumps in each interval of
time.

The key of the Harris graphical construction \cite{harris} is to show that during a certain interval of time $[0,t_0]$, $\Z^d$ can be partitioned into a countable number of finite
random sets, called islands, with no interaction between islands. For this purpose we introduce a family of random graphs containing all the information concerning the interactions needed in each interval of time $[\tau,t]$.

\paragraph{Harris random graph.} Fix $t > 0$. For each $0\leq \tau\leq t$, let $\cG_{\tau,t}=(\Z^d,\cE_{\tau,t})$ be the undirected random graph with vertex set $\Z^d$ and edge set $\cE_{\tau,t}$ defined by $\{x,y\}\in\cE_{\tau,t}$ if, and only if, $(\cT_x\cup\cT_y)\cap(\tau, t]\neq\emptyset$ and if there exists $t'\in(\cT_x\cup\cT_y)\cap (\tau,t]$ such that: (i) $y\in B(x,K_{x,n})$ if $t'=T_{x,n}$ or (ii) $x\in B(y, K_{y,m})$ if $t'=T_{y,m}$.

Note that the presence of an edge $\{x,y\}\in\cE_{\tau,t}$ indicates that a Poisson epoch has caused $x$ to look at $y$ in order to figure out how to update its spin value, or, has 
caused $y$ to look at $x$ in order to figure out how to update its spin value. Conversely, if there is no edge between $x$ and $y$ then none of them have looked at each other. The last observation implies that sites in different
components of the resulting random graph do not influence each other during the time interval $(\tau,t]$. Hence their evolutions can be computed separately.

In Section \ref{PoTP} we give the proof of the following result. 
\bteo
\label{percolation}
Let assumptions (\ref{supinf}) and (\ref{momentcondition}) hold. Then, there exists $t_0>0$ such that for any  $t\leq t_0$, the connected components of the Harris random graph $\cG_{0,t}$ are,almost surely, finite.
\eteo

Using Theorem \ref{percolation} we can show that during the time interval $[0,t_0]$, where $t_0 >0$ is deterministic and small enough, $\Z^d$ can be partitioned into a
countable number of finite islands,  with no interaction between them: $\Z^d=\cup_{\ell\in\N}C_{\ell}$, where for each $\ell\in\N, C_{\ell}$ is a finite set which is itself 
a  deterministic function of the family of Poisson processes on $[0,t_0]$ and the marks $\cK$. The collection $(C_{\ell}, \ell \in \N)$ has the additional property that 
$C_{\ell_1}\cap C_{\ell_2}=\emptyset$ for every $\ell_1\neq\ell_2$. During the time interval $[0,t_0]$ the process is constructed separately in each region $C_{\ell}$
independently of everything else. The details of this construction are given in subsection \ref{ApHGC}.

\section{Proof of Theorem \ref{coupI}}
\label{mainp}
The proof of Theorem \ref{coupI} will be divided into two steps. In the first step we introduce two families of events, $G(x,r)$ and $H(r)$, in order to study the diameter of the
cluster $C(0)$. The family of events $G(x,r)$ is helpful to understand the behavior of the cluster $C(0)$ on the subgraph of  $\cG(\cX,\cR)$ induced by the point process on $B(0,10r)$.
The family of events $H(r)$ provides a way to take care of the influence of the point process $(\cX,\cR)$ from the exterior of the ball $B(0,r)$. Our aim in this step is to show
that the probability of the percolation event can be controlled by the probabilities of the events $G(0,r)$. In the second step we will show that if the radii are not too
large, then the occurrence of the event $G(0,r_1)$ implies the occurrence of two independent events $G(x,r_2)$ and $G(x',r_2)$ where $r_1=10r_2$. Our aim in this step is to show
that the probability of the events $G(0,r_1)$ can be bounded by the square of the probability of the events $G(0,r_2)$ plus a quantity that goes to zero when $r_1$ goes to
infinity. This provides a way to take care of the probabilities of the events $G(0,r)$ that allows us to show that for $p$ small enough, $\P_{p,\,\nu}(G(0,r))$ goes to zero when
$r$ goes to infinity.

\subsection{Controlling the diameter of the cluster of the origin}
For each $x\in\Z^d$, let $D_x=\inf\{r\geq 0:\, C(x)\subset B(x,r)\}$. The percolation event is equivalent to the event $\bigcup_{x\in\Z^d}\{D_x=\infty\}$. By translation invariance of
spatially homogeneous marked point processes, the probability of the events $\{D_x>r\}$ does not depend on $x$. Therefore, the proof of Theorem \ref{coupI} is reduced to show the 
existence of $p_0>0$ such that $\lim_{r\to\infty}\P_{p,\,\nu}(D>r)=0$ for all $p<p_0$, where $D$ is the random variable $D_x$ at the origin.

We define two families of events to study the diameter of the cluster $C(0)$.   
\paragraph{The family of events $G(x,r)$.} Let $B$ be a subset of $\Z^d$. Denote by $\cG[B]$ the subgraph of $\cG(\cX,\cR)$ induced by $B$. Let $A$ be a non-empty subset of $\Z^d$
contained in $B$ and let $x\in A$. We say that $x$ is \emph{disconnected from the exterior of $A$ inside $B$} if the connected component of $\cG[B]$ containing $x$ is contained
in $A$. Now we introduce the events $G(x,r)$. Let $x\in\Z^d$ and let $r\in\N$, we say that $G(x,r)$ does not occurs if $x$ is disconnected from the exterior of $B(x,8r)$ inside $B(x,10r)$.
\paragraph{The family of events $H(r)$.} For each $r\in\N$, define
\begin{eqnarray}
\label{hr}
H(r)=\left\{\exists\,x\in \cP\cap B(0,10r)^c:\, R_x>\frac{\|x\|}{10}\right\}. 
\end{eqnarray}
The relation between the diameter of the cluster at the origin and the families of events defined above is established in the following lemma.
\blem 
\label{GHM}
The following assertion holds for all $r\in\N$:
\begin{eqnarray}
\label{Mgrande2}
G(0,r)^c\cap H(r)^c\subset\left\{D\leq 8r\right\}.
\end{eqnarray}
\elem
\paragraph{Proof of Lemma \ref{GHM}.}
If the event $H(r)$ does not occur, then there are no sites of the point process with norm greater than $10r$ connected to $B(0,9r)$. Indeed, assume that $H(r)$ does not occur. Then
for every $x\in\cP\cap B(0,10r)^c$ we have $\|x\|-R_x\geq \frac{9}{10}\|x\|>9r$. Using the triangle inequality it is easy to see that $\|y\|\geq \|x\|-R_x>9r$ for all
$y\in B(x,R_x)$. If $G(0,r)$ does not occur, then $0$ is isolated from the exterior of $B(0,8r)$. If, in addition, the event $H(r)$ does not occur, then the balls
$B(x,R_x)$ with $x\in\cP\cap B(0,10r)^c$ do not help to connect the origin to the complement of $B(0,8r)$. Thus $D\leq 8r$.  \hfill\square

\medskip
From \reff{Mgrande2} we get 
\begin{eqnarray}
\label{Mgrande3}
\P_{p,\,\nu}(D>8r)\leq \P_{p,\,\nu}(G(0,r))+\P_{p,\,\nu}(H(r)).
\end{eqnarray}
Notice that $\lim_{r\to\infty}\P_{p,\,\nu}(H(r))=0$ for all $p\in(0,1)$. This is obvious because $H(r+1)\subset H(r)$ for all $r\in\N$ and $\bigcap_{r\in\N}H(r)=\emptyset$.

\subsection{Controlling the probabilities of the events $G(0,r)$}
To take care of the probabilities $\P_{p,\,\nu}(G(0,r))$ we introduce another family of events. 
\paragraph{The family of events $\tilde{H}(r)$.} For each $r\in\N$, we define 
\begin{eqnarray}
\label{hrt}
\tilde{H}(r)=\{\exists\,x\in\cP\cap B(0,100r):\,R_x \geq r\}.
\end{eqnarray}
\blem The following inclusion holds for all $r\in\N$:
\label{Escala2}
\begin{eqnarray}
\label{escala}
G(0,10rd)\cap\tilde{H}(rd)^c&\subset&\left(\bigcup_{x\in S_{10d}}G(rx,rd)\right)\cap\left(\bigcup_{x\in S_{80d}}G(rx,rd)\right).
\end{eqnarray}
\elem
\paragraph{Proof of Lemma \ref{Escala2}.}
\label{Geom3}
Fix $r\in\N$. First, assume that the event $G(0,10rd)$ occurs but the event $\tilde{H}(rd)$ does not occur. Since $G(0,10rd)$ occurs we can go from the origin to the complement of the
ball $B(0,80rd)$ just using balls $B(x,R_x)$ centred at points from $\cP\cap B(0,100rd)$. In this way, we can go from the sphere $S_{10rd}$ to the sphere $S_{80rd}$. One of this balls,
let say $B(x_*, R_{x_*})$, touches $S_{10rd}$. Since the sphere $S_{10rd}$ is a subset of $\cup_{x \in S_{10d}}B(rx,rd)$ (see Proposition \ref{geozd} in the Appendix), we get that
this ball touches a ball of the form $B(rk, rd)$ for some $k$ in $S_{10d}$.

Now we shall prove that, for this $k$, the event $G(rk,rd)$ occurs. It is easy to see that we can go from $B(rk, rd)$ to the complement of $B(rk, 8rd)$ just using balls of the form
$B(x,R_x)$ centred at points from $\cP\cap B(0,100rd)$. Since $\tilde{H}(rd)$ does not occur, the radius of any such ball is less than $rd$. Then we can go from $B(rk,rd)$
to the complement of $B(rk,8rd)$ just using balls of the form $B(x,R_x)$ centred at points from $\cP\cap B(rk,10rd)$. In other words, the event $G(rk, rd)$occurs. Then, the event
$\bigcup_{x\in S_{10d}}G(rx,rd)$ does occur. The proof that the event $\bigcup_{x\in S_{80d}}G(rx,rd)$ does occur follows in the same lines. \hfill\square  

\medskip
The event on the right side of \reff{escala} is the intersection of two events. The first depends on what happens inside $B(0,20rd)$. The other event only depends on what happens in
the region $B(0,70rd)^c$. Then, these two events are independent. By translation invariance of spatially homogeneous marked point processes we get
\begin{eqnarray}
\label{Go8bis}
\P_{p,\,\nu}(G(0,10rd))\leq |S_{10d}||S_{80d}|\P_{p,\,\nu}(G(0,rd))^2+\P_{p,\,\nu}(\tilde{H}(rd)).
\end{eqnarray}

\blem
\label{Cotas} 
There exist positive constants $C_2$ and $C_3$, which depends only on the dimension $d$, such that, for any $r\in\N$, the following inequalities hold:
\begin{eqnarray}
\label{CF1b}
\P_{p,\,\nu}(G(0,r)) &\leq& p\, C_2r^d,\\
\label{CG1b}
\P_{p,\,\nu}(\tilde{H}(r))&\leq& p\,C_3\E_{p,\,\nu}\left[R^d\one\{R\geq r\}\right].
\end{eqnarray}
\elem
\paragraph{Proof of Lemma \ref{Cotas}.}
It is a simple geometric fact that there exists a positive constant $C$ which depends only on the dimension $d$ such that $|B(0,r)|\leq Cr^d$.
 
Let $r\in\N$. A simple computation shows that
\begin{eqnarray}
\P_{p,\,\nu}(G(0,r))
&\leq& \P_{p,\,\nu}(\exists\,x\in\cP\cap B(0,10r))\nonumber\\
&\leq& p\,|B(0,10r)|.
\end{eqnarray}
The inequality \reff{CF1b} is satisfied with $C_2=10^dC$. 

To show \reff{CG1b} we note that $\tilde{H}(r)=\one\{X\geq 1\}$, where $X$ is a random variable defined by
\[X=\sum_{x\in B(0,100r)}\one\{x\in\cP\}\one\{R_x\geq r\}.\]
We have 
\begin{eqnarray*}
\P_{p,\,\nu}(\tilde{H}(r))&\leq& \E_{p,\,\nu}\left[X\right] \nonumber \\
&=&\sum_{x\in B(0,100r)}p\,\P_{p,\,\nu}(R_x\geq r)\nonumber \\
&=& p\, \left|B(0,100r)\right|\P_{p,\,\nu}(R\geq r)\nonumber \\
&\leq& p\,C_3\E_{p,\,\nu}[R^d\one\{R\geq r\}],
\end{eqnarray*}
where $C_3=100^dC$. The first equality follows from the independence between $\cP$ and $\cR$ and the second equality follows from the fact that the random variables $(R_x, x\in\Z^d)$
are identically distributed. \hfill\square

\subsection{Proof of Theorem 1}
By \reff{Mgrande3}, the proof of Theorem \ref{coupI} is reduced to show the existence of $p_0>0$ such that there exists an increasing sequence $(r_n)_{n\in\N}\subset \N$ with
$\lim_{n\to\infty}\P_{p,\,\nu}(G(0,r_n))=0$ for any $p<p_0$. For this reason we need the following lemma.
\blem 
\label{FG0}
Let $f$ and $g$ be two functions from $\N$ to $\R_+$ satisfying the following conditions: (i) $f(r)\leq 1/2$ for all $r\in \{1,\dots, 10\}$; (ii) $g(r)\leq 1/4$ for all $r\in\N$; (iii) for all $r\in\N$:     
\begin{eqnarray}
\label{FG1}
f(10r)\leq f^2(r)+g(r).
\end{eqnarray}
If $\lim_{r\to\infty}g(r)=0$, then $\lim_{n\to\infty}f(10^nr)=0$ for each $r\in\{1,\dots, 10\}$. 
\elem

\paragraph{Proof of Lemma \ref{FG0}.} For each $n\in\N$, let  $F_n=\max_{1\leq r\leq 10}f(10^nr)$ and let $G_n=\max_{1\leq r \leq 10}g(10^nr)$. Using \reff{FG1} and hypothesis (i) and (ii) we may conclude, by means of the induction principle that, for each $n\in\N$,  $F_n\leq 1/2$ and 
\begin{eqnarray}
\label{CotaIndb}
F_n \leq \frac{1}{2^{n+1}}+\displaystyle \sum_{j=0}^{n-1}\frac{1}{2^j}G_{n-1-j}. 
\end{eqnarray}
Since $g(10^nr)$ goes to zero as $n\to\infty$ we have that $G_n\to 0$ when $n\to\infty$. By \reff{CotaIndb}, we obtain that $F_n\to 0$ when $n\to\infty$. \hfill\square

\medskip

Consider the functions $f(r)=C_1\P_{p,\,\nu}(G(0,rd))$ and $g(r)=C_1\P_{p,\,\nu}(\tilde{H}(rd))$, where $C_1=|S_{10d}||S_{80d}|$. By \reff{Go8bis}, it follows that 
\begin{eqnarray}
\label{Go8tis}
f(10r)\leq f^2(r)+g(r).
\end{eqnarray}
By condition $\E_{p,\nu}[R^d]=\sum_{r\geq 1}r^d\nu(r)<\infty$ and \reff{CG1b}, we have that $\lim_{r\to\infty}g(r)=0$ for any $p$.

We show that there exists $p_0>0$ such that if $p<p_0$ then $f(r)\leq 1/2$, $1\leq r\leq 10$ and $g(r)\leq 1/4$, $r\in\N$. 

Set
\[
p_0=\min((2C_1C_2(10d)^d)^{-1}, (4C_1C_3\E_{p,\,\nu}[R^d])^{-1}).
\]

By condition $\E_{p,\,\nu}[R^d]<\infty$, we get $p_0>0$.

Let $p>0$ be such that $p\leq p_0$. It follows from\reff{CF1b} that 
\[
f(r)\leq \frac{1}{2}\left(\frac{r}{10}\right)^d.
\]
Thus we have that if $0<p\leq p_0$, then $\max_{1\leq r\leq 10}f(r)\leq 1/2$. 

By \reff{CG1b}, we get
\[
g(r)\leq \frac{1}{4}.
\]

Finally, by Lemma \ref{FG0}, we may conclude that $\lim_{n\to\infty}f(10^nr)=0$ for each $r\in\{1,\dots, 10\}$. In particular,
\[\lim_{n\to\infty}f(10^n)=\lim_{n\to\infty}C_1\P_{p,\,\nu}(G(0, 10^nd)=0.\]
\hfill\square

\bigskip

We finish this section by proving the complete coverage of $\Z^d$ under the assumption $\E_{p,\,\nu}[R^d]=\infty$.

\subsection{Proof of Theorem \ref{Riidinf}}
We prove the equivalent statement that, for all $r\in\N$, the following assertion holds:
\[\P_{p,\,\nu}(\exists\, x\in\cP: B(0,r)\subset B(x,R_x))=1.\]
If $R_x>\|x\|+r$, then $B(0,r)\subset B(x,R_x)$. Hence,
\begin{eqnarray*}
\P_{p,\,\nu}(\exists\, x\in\cP: B(0,r)\subset B(x,R_x))\geq\P_{p,\,\nu}(\exists\, x\in\cP: R_x>\|x\|+r) 
\end{eqnarray*}
Let $A_k$ be the event defined by $A_k=\{\exists x\in\cP\cap S_k: \, R_x>k+r\}$. It is clear that the events  $A_k$ are independent and that $\P_{p,\,\nu}(A_k)=p|S_k|\P_{p,\,\nu}(R>k+r)$.
Note that
\begin{eqnarray}
\label{ineq}
\sum_{k\geq 0}\P_{p,\,\nu}(A_k)&=&p\sum_{k\geq 0}|S_k|\P_{p,\,\nu}(R>k+r)\nonumber \\
&=&\sum_{k\geq 0}|B_k|\P_{p,\,\nu}(R=k+r+1).
\end{eqnarray}
Since $\E_{p,\,\nu}[R^d]=\infty$, we conclude that the series in the right hand side of \reff{ineq} diverges. By the second Borel-Cantelli lemma (see Durrett \cite{durrett2}, page 50),
we have that $\P_{p,\,\nu}(A_k$ i.o.$)=1$. \hfill\square 

\section{Proof of Theorem \ref{percolation}}
\label{PoTP}

The proof of Theorem \ref{percolation} falls naturally into two steps. In the first step, we construct a family of marked point processes $(\cX_t,\cR_t)$ on $\Z^d$ such that the random graphs $\cG(\cX_t,\cR_t)$ and $\cG_{0,t}$ have the same distribution. In the second step, using Theorem \ref{coupI}, we show that for $t$ small enough the connected components of the random graph $\cG(\cX_t,\cR_t)$ are, almost surely, finite. For the sake of clarity, each step is divided into a sequence of lemmas.
\bigskip

In order to prove Theorem \ref{percolation}, we need to introduce some notation. For each $x\in\Z^d$, $r\in \{-1,0,1,2,\dots\}$ and $0< t$ let
\[
N_{x,r}(t):=\sum_{n\geq 1}\one\{T_{x,n}\leq t\}\one\{K_{x,n}>r\}.
\]

$N_{x,r}(t)$ is nothing but the number of occurrences of the marked Poisson process $(\cT_x, \cK_x)$ during the time interval $(0,t]$ whose marks are greater than $r$. Notice that 
$N_{x}(t):=N_{x,-1}(t)$ is the counting measure associated to the Poisson process $\cT_x$. 

Let $(\cT,\cK)=\{(\cT_x, \cK_x):\,x\in\Z^d\}$ be a family of mutually independent marked Poisson point processes on the time line $[0,\infty)$. Let $\P_{(\cT,\cK)}$ and 
$\E_{(\cT,\cK)}$ respectively be the probability measure and the expectation operator induced by $(\cT,\cK)$.

\begin{rem}
For each $x\in\Z^d$, let $M_x(t):=\sum_{n\geq 1}\max(K_{x,1},\dots, K_{x,n})\one\{N_x(t)=n\}$. It follows from the construction described above and the Coloring Theorem (see Kingman
\cite{Kingman}, page 52) that
\begin{eqnarray}
\label{M1}
\P_{(\cT,\cK)}(M_x(t)\leq r)=\P_{(\cT,\cK)}(N_{x,r}(t)=0)=\exp\left(-M_xtG_x(r)\right),
\end{eqnarray}
where $G_x(r)=\sum_{k>r}\nu_{x}(k)$.
Also, we have
\begin{eqnarray}
\label{M2}
\P_{(\cT,\cK)}(M_x(t)\leq r|N_{x}(t)\geq 1)=\frac{\exp\left(-M_xtG_{x}(r)\right)-\exp\left(-M_xt\right)}{1-\exp\left(-M_xt\right)}.
\end{eqnarray}
\end{rem}

\subsection{Simultaneous coupling construction: Step 1.}
Let $\cU_1=(U_{x,1}:\,x\in\Z^d)$ and $\cU_2=(U_{x,2}:\,x\in \Z^d)$ be two mutually independent families of independent uniform random variables on $(0,1]$. The probability space where
this coupling is performed is the one where the two families of uniform random variables are defined. We denote by $\P_{(\cU_1,\cU_2)}$ the probability measure induced by the families
of random variables $\cU_1$ and $\cU_2$.

For each $x\in\Z^d$ and $t > 0$, define
\begin{eqnarray}
\label{PP1}
X_{x,t}&:=&\one\{U_{x,1}\leq 1-\exp(-M_xt)\},\\
\label{MP1}
R_{x,t}&:=&F^{-1}_{x,t}(U_{x,2}),
\end{eqnarray}
where $F_{x,t}(r)$ is the cumulative distribution function given by
\begin{eqnarray}
\label{DD1bis}
F_{x,t}(r)=\frac{\exp\left(-M_xtG_{x}(r)\right)-\exp\left(-M_xt\right)}{1-\exp\left(-M_xt\right)}.
\end{eqnarray}

Set $\cX_t=(X_{x,t}:\,x\in\Z^d)$ and $\cR_t=(R_{x,t}:\,x\in\Z^d)$. It follows from the construction that the process $(\cX_t,\cR_t)$ is a marked point process on $\Z^d$ satisfying
\begin{eqnarray}
\label{DPbis}
\P_{(\cU_1,\cU_2)}(X_{x,t}=1)&=&\P_{(\cT,\cK)}(N_x(t)\geq 1),\\
\label{DD1b}
\P_{(\cU_1,\cU_2)}(R_{x,t}\leq r)&=&\P_{(\cT,\cK)}(M_x(t)\leq r|N_{x}(t)\geq 1).
\end{eqnarray}

\blem
\label{Lpaso1} Let $t > 0$ and let $(\cX_t,\cR_t)$ be the marked point process on $\Z^d$ defined by \reff{PP1} and \reff{MP1}. Then, the random graphs $\cG(\cX_t,\cR_t)$ and $\cG_{0,t}$ are equally distributed.
\elem

\paragraph{Proof of Lemma \ref{Lpaso1}.}
By \reff{DPbis} and \reff{DD1b} we have that 
\begin{eqnarray}
\label{Leyes1b}
\P_{(\cU_1,\cU_2)}(X_{x,t}=1, R_{x,t}\leq r)=\P_{(\cT,\cK)}(N_x(t)\geq 1, M_x(t)\leq r). 
\end{eqnarray}
Then, the random graphs $\cG_{0,t}(\cT,\cK)$ and $\cG(\cX_t,\cR_t)$ have the same distribution. \hfill\square

\subsection{Properties of the coupled processes $(\cX_t,\cR_t)$}
Now, we study the main properties of the marked point processes $(\cX_t,\cR_t)$ needed for the proof of Theorem \ref{percolation}. We begin by proving the following auxiliary result.

\blem
\label{haux}
For any $a\in(0,1)$, the function
\begin{eqnarray}
h(z)=\frac{1-\exp(-az)}{1-\exp(-z)}
\end{eqnarray}
is non-decreasing on $[0,\infty)$.
\elem

\paragraph{Proof of Lemma \ref{haux}.}
It suffices to prove the result for rational $a$. Then, assume that $a$ is a positive rational number and write it as the ratio of two positive integers $a=m/n$, where $1<m<n$. Now,
making $y=\exp(-z/n)$ we get
\begin{eqnarray}
\frac{1-\exp(-az)}{1-\exp(-z)}&=&\frac{1-y^m}{1-y^n}=\left(1+\frac{\sum_{i=0}^{n-m-1}y^k}{\sum_{k=1}^my^{-k}}\right)^{-1}.
\end{eqnarray}
Since the expression above is a decreasing function of $y$ and $y=\exp(-z/n)$ is itself a decreasing function of $z$, the result follows.  \hfill\square     

\blem
\label{SDom}
If $0<t'< t\leq 1$, then
\begin{eqnarray}
(\cX_{t'}, \cR_{t'})\preceq (\cX_t, \cR_t).
\end{eqnarray}
\elem

\paragraph{Proof of Lemma \ref{SDom}.}
Fix $x\in\Z^d$. It follows from \reff{PP1} that $X_{x,t'}\leq X_{x,t}$. Fix $r\in\N_0$. We shall show that if $0<t'<t$, then
\begin{eqnarray}
\label{losradioscrecen}
F_{x,t}(r)\leq F_{x,t'}(r).
\end{eqnarray}
By \reff{DD1bis}, it suffices to show that 
\begin{eqnarray}
\tilde{h}(t)=\frac{1-\exp\left(-M_xtG_{x}(r)\right)}{1-\exp\left(-M_xt\right)}
\end{eqnarray}
is a non-decreasing function on $(0,1]$. This follows from Lemma \ref{haux} by substituting $a$ by $G_x(r)$ and $z$ by $M_xt$.

From \reff{losradioscrecen} we may conclude that the random variables $R_{x,t}$ defined in \reff{MP1} satisfy  $R_{x,t'}\leq R_{x,t}$. hfill\square

\subsection{Simultaneous coupling construction: Step 2.}
For each $x\in\Z^d$ and $t > 0$, let
\begin{eqnarray}
\label{PP2}
X^*_{x,t}&:=&\one\{U_{x,1}\leq 1-\exp(-M^*t)\},\\
\label{hmp}
R^*_{x}&:=&\sum_{r\in\mathbb{N}_0}r\one\left\{\inf_{x\in\Z^d}F_{x,1}(r-1)<U_{x,2}\leq \inf_{x\in\Z^d}F_{x,1}(r) \right\}
\end{eqnarray}
and define $\cX^*_t=(X^*_{x,t}:\,x\in\Z^d)$ and $\cR^*=(R^*_{x}:\,x\in\Z^d)$. 

First, we use conditions \reff{supinf} and \reff{momentcondition} to show that the random variables introduced in \reff{hmp} are well defined. For that purpose it suffices to show that $\lim_{r\to\infty}\inf_{x\in\Z^d}F_{x,1}(r)=1$. By \reff{supinf}, we have
\begin{eqnarray*}
F_{x, 1}(r)&=&\frac{\exp\left(-M_xG_{x}(r)\right)-\exp\left(-M_x\right)}{1-\exp\left(-M_x\right)}\nonumber\\
&=&1-\frac{1-\exp\left(-M_xG_{x}(r)\right)}{1-\exp\left(-M_x\right)}\nonumber\\
&\geq&1-\frac{1-\exp\left(-M^*\sup_{x\in\Z^d}G_{x}(r)\right)}{1-\exp\left(-M_*\right)}
\end{eqnarray*}
for any $x\in\Z^d$. Then,
\begin{eqnarray}
\label{CIFx1}
\inf_{x\in\Z^d}F_{x, 1}(r)&\geq&1-\frac{1-\exp\left(-M^*\sup_{x\in\Z^d}G_{x}(r)\right)}{1-\exp\left(-M_*\right)}.
\end{eqnarray}
On the other hand, 
\begin{eqnarray}
\label{CSup}
\sup_{x\in\Z^d}G_{x}(r)=\sup_{x\in\Z^d}\sum_{\ell>r}\nu_x(\ell)\leq\sum_{\ell>r}\sup_{x\in\Z^d}\nu_x(\ell).
\end{eqnarray}
Under condition \reff{momentcondition} we may conclude that the right hand side of  \reff{CSup} converges to $0$ when $r$ goes to $\infty$. Then
\begin{eqnarray}
\label{CSup2}
\lim_{r\to\infty}\sup_{x\in\Z^d}G_{x}(r)=0. 
\end{eqnarray}
From \reff{CIFx1} and \reff{CSup2} we deduce that $\lim_{r\to\infty}\inf_{x\in\Z^d}F_{x,1}(r)=1$.

\subsection{The dominating marked point process $(\cX^*_t, \cR^*)$}
\blem
\label{DMPP}
For any $0<t\leq 1$,
\begin{eqnarray}
\label{DSt1}
(\cX_t, \cR_t) \preceq (\cX^*_t, \cR^*).
\end{eqnarray}
\elem

\paragraph{Proof of Lemma \ref{DMPP}.}
It follows from the construction that $X_{x,t}\leq X^*_{x,t}$ and $R_{x,1}\leq \R^*_x$. By Lemma \ref{SDom}, we have that $R_{x,t}\leq R_{x,1}$. This completes the proof. 
\hfill\square

The last ingredient needed to prove Theorem \ref{percolation} is the following result.

\blem 
\label{Cota1} Under the same assumptions of Theorem \ref{percolation} there exists $0<t_0\leq 1$ such that $\P_{(\cU_1,\cU_2)}(Percolation)=0$ for all $0<t\leq t_0$. 
\elem
\paragraph{Proof of Lemma \ref{Cota1}.}
First note that, for each $t>0$, $(\cX^*_t, \cR^*)$ is a spatially homogeneous, marked point processes on $\Z^d$ with retention parameter $p(t)=1-\exp(-M^*t)$ and probability function of its marks $\nu(r)=\P_{(\cU_1,\cU_2)}(\hR_x=r)$ satisfying
\begin{eqnarray}
\label{Cotagrosa}
\sum_{r\in\N_0}r^d\nu(r)<\infty.
\end{eqnarray}
Indeed,  
\begin{eqnarray}
\label{Cotagrosa1}
\nu(r)&=&\inf_{x\in\Z^d}\P_{(\cU_1,\cU_2)}(R_{x,1}\leq r)-\inf_{x\in\Z^d}\P_{(\cU_1,\cU_2)}(R_{x,1}\leq r-1)\nonumber\\
&\leq& \sup_{x\in\Z^d}\P_{(\cU_1,\cU_2)}(R_{x,1}=r)
\end{eqnarray}
Inequality \reff{Cotagrosa1} follows from the inequality $\inf\{a_x+b_x\}\leq \inf\{a_x\}+\sup\{b_x\}$  applied to the sequence $a_{x}(r)=\P_{(\cU_1,\cU_2)}(R_{x,1}\leq r-1)$ and $b_x(r)=\P_{(\cU_1,\cU_2)}(R_{x,1}=r)$.

On the other hand, we have 
\begin{eqnarray}
\label{Cotagrosa2}
\P_{(\cU_1,\cU_2)}(R_{x,1}=r)
&=&\P_{(\cU_1,\cU_2)}(R_{x,1}\leq r)-\P_{(\cU_1,\cU_2)}(R_{x,1}\leq r-1)\nonumber\\
&=&\frac{\exp\left(-M_xG_{x}(r)\right)-\exp\left(-M_xG_{x}(r-1)\right)}{1-\exp\left(-M_x\right)}\nonumber\\
&=&\frac{\exp\left(-M_xG_{x}(r)\right)}{1-\exp\left(-M_x\right)}\left(1-\exp\left(-M_x\nu_{x}(r)\right)\right)\nonumber\\
&\leq&\left(\frac{M^*}{1-\exp\left(-M_*\right)}\right)\nu_{x}(r).
\end{eqnarray}
The last inequality follows from well known properties of the exponential function.

From \reff{Cotagrosa1} and \reff{Cotagrosa2} we may conclude that 
\begin{eqnarray}
\label{Cotagrosa3}
\nu(r)\leq \left(\frac{M^*}{1-\exp\left(-M_*\right)}\right)\sup_{x\in\Z^d}\nu_{x}(r).
\end{eqnarray}

Inequality \reff{Cotagrosa} follows immediately from \reff{Cotagrosa3}.

Since $(\cX_t^*, \cR^*)$ satisfies the hypothesis of Theorem \ref{coupI}, there exists $p_0>0$ such that $\P_{(\cU_1,\cU_2)}($Percolation$)=0$ for any $0<t\leq 1$ such that $p(t)\leq p_0$.
Therefore, the connected components of the random graphs $\cG(\cX^*_t, \cR^*)$ are almost surely finite for any $0<t\leq \min(-\frac{1}{M^*}\log(1-p_0), 1)$. Indeed,
\begin{eqnarray*}
1-\exp(-M^*t)\leq p_0&\iff& t\leq -\frac{1}{M^*}\log(1-p_0). 
\end{eqnarray*}
\hfill\square

\subsection{Proof of Theorem \ref{percolation}.} 
By Lemma \ref{Lpaso1}, the random graphs $\cG_{0,t}(\cT,\cK)$ and $\cG(\cX_t,\cR_t)$ have the same distribution. By Lemmas \ref{DMPP} and \ref{Cota1}, we have that there
exists $0<t_0\leq 1$ such that $\P_{(\cU_1,\cU_2)}($Percolation$)=0$ for all $0<t\leq t_0$. Therefore, the connected components of the Harris random graph
$\cG_{0,t}(\cT,\cK)$ are almost surely finite for all $0<t\leq t_0$. \hfill\square

\section{Appendix}

\subsection{Harris graphical construction}
\label{ApHGC}
Let $t_0>0$ be as in Theorem \ref{percolation} and let $C_{\ell}=C_{\ell}(\cT\cap[0,t_0],\cK)$, $\ell\in\N$, be the partition of $\Z^d$ into a countable number of finite islands, with no interaction between them.   

\paragraph{Finite-volume construction}
The construction of the process on each finite island $C_{\ell}$ with initial configuration $\eta$ using the Poisson processes is straightforward because the epoch of the associated Poisson processes are well ordered.  

Fix $\ell\in\N$ and consider $0<\tau_1<\tau_2<\cdots<\tau_{n}$, where
\begin{eqnarray}
\{\tau_1,\tau_2, \dots, \tau_{n}\}=\bigcup_{x\in C_{\ell}}\left(\cT_x\cap [0,t_0]\right).
\end{eqnarray}
We construct the process $\sigma^{\eta}_t(C_\ell)$ inductively as follows. For $k=1,\dots, n$, let $x_1,x_2,\dots,x_n$ be the sites such that $\tau_k\in\cT_{x_k}$ .
\paragraph{Step 1.} Let
\[\sigma^{\eta}_t(C_\ell)=\eta(C_\ell) \mbox{ for all } 0\leq t<\tau_1\]
and set
\begin{eqnarray*}
\sigma^{\eta}_{\tau_1}(C_\ell)=\sigma^{\eta}_{x_1, W_{x_1}(\sigma^{\eta}_{\tau_1-}(C_\ell))}.
\end{eqnarray*}
Thus we have defined the process on the time interval $[0,\tau_1]$.
\paragraph{Inductive step.} Assume that $\sigma^{\eta}_{t}(C_\ell)$ has already been defined for all $0\leq t\leq \tau_k$. Then set
\[\sigma^{\eta}_{t}(C_\ell)=\sigma^{\eta}_{\tau_k}(C_\ell) \mbox{ for all } \tau_k<t<\tau_{k+1}\]
and
\begin{eqnarray*}
\sigma^{\eta}_{\tau_{k+1}}(C_\ell)=\sigma^{\eta}_{x_{k+1}, W_{x_{k+1}}(\sigma^{\eta}_{\tau_{k+1-}}(C_\ell))}.
\end{eqnarray*}

This step is repeated until the construction has been finished on the time interval $[0,t_0]$.

\paragraph{Infinite-volume construction.} Finally, for each $\ell\in\N$, let 
\begin{eqnarray}
\label{GSb}
(\sigma^\eta_t)(C_{\ell}):=\sigma^{\eta}_t(C_\ell), \qquad t\in(0,t_0].
\end{eqnarray}
Since $t_0$ is independent of the initial configuration, we may conclude, by means of the Markovian property of Poisson processes,  that the state of the process may be computed by
induction at any time $t \geq 0$. 

\subsection{Geometry of $\Z^d$}
The following proposition deals with geometric aspects used in the proof of Lemma \ref{Escala2}.
\bpro
\label{geozd}
Fix $d\in\N$. Then, for any $n,r\in\N$, we have
\[S_{nr}\subset\bigcup_{x \in S_n} B\left(rx,\frac{d}{2}\,r\right).\]
\epro

In order to prove Proposition \ref{geozd} we need the following result.

\blem 
\label{geozd2}
Let $x=(x_1,\dots, x_d)\in\R^d$ be such that $1>x_1\geq x_2\geq\cdots\geq x_d>0$ and $\sum_{i=1}^dx_i=m\in\N$, where $m<d$. Let $y=(1,\dots, 1, 0,\dots, 0)$ ($m$-ones). Then, $\|y-x\|\leq \frac{d}{2}$.
\elem

\paragraph{Proof of Proposition \ref{geozd}.} By symmetry, it suffices to prove the proposition for $x$ in the region
$\{x=(x_1, x_2, \ldots, x_d) \in \mathbb{R}^d : x_i \geq0, i=1, 2, \ldots, d \}$. Fix $n$ and $r\in\N$. By Lemma \ref{geozd2}, we have that for any $x\in\R^d$ with
$\|x\|=n$ there exists $y\in\Z^d$ with $\|y\|=n$ such that $\|y-x\|\leq\frac{d}{2}$. 

Pick $x\in S_{nr}$. Then, for $\frac{x}{r}\in S_n$ there exists $ y\in\Z^d$ with $\|y\|=n$ such that $\|y-\frac{x}{r}\|\leq \frac{d}{2}$. Thus, $\|ry-x\|\leq\frac{d}{2}r$.
\hfill\square
\vspace{0.5cm}

We finish this subsection by proving Lemma \ref{geozd2}. 
\paragraph{Proof.} We begin by observing that
\begin{eqnarray}
\label{n162}
\|y-x\|
&=&\sum_{i=1}^m(1-x_i)+\sum_{i=m+1}^dx_i\nonumber\\
&=&1-x_1+\sum_{i=2}^m(1-x_i)+\sum_{i=m+1}^dx_i\nonumber\\
&=&1-m+\sum_{i=2}^dx_i+m-1-\sum_{i=2}^mx_i+\sum_{i=m+1}^dx_i\nonumber\\
&=&2\sum_{i=m+1}^dx_i.
\end{eqnarray}
Now assume that $\|y-x\|>\frac{d}{2}$. From \reff{n162} we have that $\sum_{i=m+1}^dx_i>\frac{d}{4}$. Therefore, there exists $i\in\{m+1,\dots, d\}$ such that $x_i>\frac{d}{4(d-m)}$.
Also, $x_1,\dots, x_m>\frac{d}{4(d-m)}$. Since $\sum_{i=1}^mx_i>\frac{dm}{4(d-m)}$ we get
\begin{eqnarray}
\label{n163}
\sum_{i=1}^dx_i>\frac{dm}{4(d-m)}+\frac{d}{4}
=\frac{d^2}{4(d-m)}.
\end{eqnarray}
Note that 
\begin{eqnarray}
\label{n164}
\frac{d^2}{4(d-m)}\geq m
\end{eqnarray}
if, and only if, $(d-2m)^2\geq 0$. Since the last inequality is true, so is inequality (\ref{n164}). From \reff{n163} and \reff{n164} we get that $\sum_{i=1}^dx_i>m$ which is a contradiction. 
The contradiction comes from the fact that we have assumed that $\|y-x\|>\frac{d}{2}$. \hfill\square 

\section*{Acknowledgement}
We thank Jorge R. Busch for a carefully reading of a previous version of this work and for comments that improved the presentation of the results. During the realization of this work
both authors received partial financial support from FAPESP, grant 09/52379-8. Also, the second author received support from FAPESP, grant  2009/16437-3.

\end{document}